\documentclass[11pt]{article}
\usepackage[T1]{fontenc}
\usepackage{amsmath,amssymb,amsthm,mathtools}
\usepackage{lmodern,microtype}
\usepackage[margin=1in]{geometry}
\usepackage[hidelinks]{hyperref}
\numberwithin{equation}{section}
\newtheorem{theorem}{Theorem}[section]
\newtheorem{lemma}[theorem]{Lemma}
\newtheorem{corollary}[theorem]{Corollary}
\newtheorem{proposition}[theorem]{Proposition}
\theoremstyle{definition}

\theoremstyle{remark}

\newcommand{\N}{\mathbb N}
\newcommand{\R}{\mathbb R}
\newcommand{\Tr}{\mathrm{Tr}}
\newcommand{\IF}{\mathrm{IF}}
\newcommand{\WF}{\mathrm{WF}}
\newcommand{\id}{\mathrm{id}}
\newcommand{\rank}{\operatorname{rank}}
\newcommand{\graph}{\operatorname{graph}}
\newcommand{\sym}{\operatorname{sym}}
\newcommand{\Bred}{\leq_B}
\title{No Borel upper bound for the $F_\sigma$ equivalence relations}
\author{Ruijun Wang}
\date{}
\begin{document}
\maketitle
\begin{abstract}
We prove that the family of coset relations of $F_\sigma$ real linear subspaces of Banach space $\ell_1$ has no common upper bound among Borel equivalence relations in the sense of Borel reducibility. These relations are indexed by a well-founded tree. Consequently, on a Polish space, a universal \(\Sigma^0_\alpha\) equivalence relation exists exactly when \(\alpha=1\), and a universal \(\Pi^0_\alpha\) equivalence relation exists exactly when \(\alpha=1\) or \(2\). This strengthens the earlier nonexistence theorem for a universal \(F_\sigma\) equivalence relation. This also answers a question of Kanovei and Lyubetsky about universal equivalence relations in Borel classes.
\end{abstract}

\section{Introduction}\label{sec:introduction}

A classification problem asks for invariants which determine when two
objects are equivalent. Invariant descriptive set theory makes the
regularity of such an assignment part of the problem. The objects are
parametrized by standard Borel spaces, and a Borel reduction is a
Borel assignment of invariants which preserves and reflects equivalence.
Thus \(E\Bred F\) means that classification up to \(E\) can be carried
out using invariants classified up to \(F\). This viewpoint connects
classification questions in logic, topology, algebra, and analysis.

The theory grew from classical descriptive set theory and the study of
orbit equivalence relations. Friedman and Stanley's work on countable
structures \cite{FriedmanStanley} established Borel reducibility as a
systematic framework for comparing classification problems. The
Harrington--Kechris--Louveau dichotomy \cite{HKL} identified eventual
equality as an obstruction to classification by real numbers. Gao's
monograph \cite{Gao} develops these methods, Polish group actions,
classification by countable structures, and the benchmark equivalence
relations used to compare concrete problems. Kechris's book \cite{Kechris}
provides the underlying descriptive set theory.

One basic question asks how the Borel complexity of an equivalence
relation constrains its position under Borel reducibility. For a Borel
pointclass \(\Gamma\), a \emph{universal} (or \emph{complete})
\(\Gamma\) equivalence relation must itself belong to \(\Gamma\) and
must dominate every equivalence relation in that class. This is a
question about reductions induced by maps of the underlying spaces;
it is different from completeness of a set under arbitrary continuous
preimages.

Rosendal \cite[Definition~18 and Proposition~19]{Rosendal} exhibited a
universal \(K_\sigma\) degree, represented by bounded coordinatewise
difference on a compact product of finite sets and Borel bireducible
with \(\ell_\infty\). He also asked whether every \(F_\sigma\)
equivalence relation reduces to that degree
\cite[Problem~28]{Rosendal}. Here \(K_\sigma\) means a countable union
of compact sets, whereas \(F_\sigma\) means a countable union of closed
sets. On a compact space the two conditions coincide, but they need
not coincide on an arbitrary Polish space.

The earlier paper \cite[Theorem~1.1]{WangFsigma} proves that for every
\(F_\sigma\) equivalence relation \(F\), there is an \(F_\sigma\)
real linear subspace \(H\subseteq\ell_2\) whose coset relation has no
Baire-measurable reduction to \(F\). In particular, no \(F_\sigma\)
equivalence relation is universal for the \(F_\sigma\) equivalence
relations on arbitrary Polish spaces. That theorem leaves open the
possibility of a common upper bound of higher Borel complexity.

Kanovei and Lyubetsky \cite[Section~9, p.\ 102]{KanoveiLyubetsky}
asked which Borel classes contain complete equivalence relations.
Their discussion of Borel-class estimates concerns Cantor-type spaces.
They construct a cofinal family of Borel equivalence relations, but
cofinality across the whole hierarchy does not by itself give a
universal relation at any prescribed level. We address both the
arbitrary-Polish-space question and their higher-level Cantor-space
formulation.

Our main result excludes an upper bound of any Borel complexity.

\begin{theorem}\label{thm:main}
There is no Borel equivalence relation \(F\) such that every
\(F_\sigma\) equivalence relation on a Polish space Borel reduces to \(F\).
\end{theorem}

In fact, the counterexamples can all be taken to be coset relations
of \(F_\sigma\) real linear subspaces of one fixed space isomorphic
to \(\ell_1\); see Corollary~\ref{cor:linear}. The improvement over
\cite{WangFsigma} is the unrestricted Borel complexity of the target.
The proof here is independent of the proof of that earlier theorem.

Since every \(F_\sigma\) set belongs to \(\Sigma^0_\alpha\) for
\(\alpha\geq2\) and to \(\Pi^0_\alpha\) for \(\alpha\geq3\),
Theorem~\ref{thm:main} gives the negative entries in the following
complete classification. The positive entries are classical.
\[
\begin{array}{c|c|c}
\alpha & \text{universal }\Sigma^0_\alpha
       & \text{universal }\Pi^0_\alpha\\ \hline
1 & {=}_{\N} & {=}_{2^{\N}}\\
2 & \text{none} & {=}_{2^{\N}}\\
3\leq\alpha<\omega_1 & \text{none} & \text{none}
\end{array}
\]
In this table domains range over all Polish spaces. For every
\(\alpha\geq3\), both negative conclusions remain true when all
domains are required to be Cantor space; see
Corollary~\ref{cor:cantor}. Compact \(F_\sigma\) relations retain the
universal \(K_\sigma\) degree, so that distinction matters at level two.

This paper is a collaboration with artificial intelligence. The AI model is gpt 6. The AI proves the main theorem under human suggestions and writes the main body of early drafts. The author revises the paper and improves readability. The author has reviewed all AI-generated content.

The rest of the paper is organized as follows. Section~\ref{sec:preliminaries} gives the definitions and classical
facts needed in the proof. Section~\ref{sec:transfer} proves a category
transfer lemma: a common Borel bound would force a certain analytic
membership test to be Borel. Section~\ref{sec:tensor} establishes an
\(F_\sigma\) image lemma for symmetric finite-rank tensors.
Section~\ref{sec:trees} combines it with tree kernels and a fixed
two-step nilpotent Polish group to prove Theorem~\ref{thm:main}.
Section~\ref{sec:linear} passes to additive linear cosets, and
Section~\ref{sec:hierarchy} derives the Borel-hierarchy conclusions.

\section{Preliminaries}\label{sec:preliminaries}

\subsection{Polish spaces and Borel classes}
A \emph{Polish space} is a separable completely metrizable topological
space. Its Borel sets form the smallest \(\sigma\)-algebra containing
the open sets. A \emph{standard Borel space} is a set equipped with
the Borel \(\sigma\)-algebra of some Polish topology. A map is Borel
if preimages of open sets are Borel. We use
\(\N=\{0,1,2,\ldots\}\), Cantor space \(2^{\N}\), and Baire space
\(\N^{\N}\), with their product topologies.

For a Polish space \(X\), let \(\Sigma^0_1(X)\) be the open sets and
\(\Pi^0_1(X)\) the closed sets. For \(1<\alpha<\omega_1\), define
\(\Sigma^0_\alpha(X)\) to consist of countable unions of sets in
\(\bigcup_{\beta<\alpha}\Pi^0_\beta(X)\), where \(\beta\geq1\), and
let \(\Pi^0_\alpha(X)\) consist of their complements. Write
\(\Delta^0_\alpha=\Sigma^0_\alpha\cap\Pi^0_\alpha\).
In particular, \(\Sigma^0_2\) is the class of \(F_\sigma\) sets and
\(\Pi^0_2\) the class of \(G_\delta\) sets. These classes are
preserved by continuous preimages. Their standard closure and
relativization properties will be used below \cite{Kechris}.
The assertion that an equivalence relation on \(X\) belongs to one
of these classes refers to its graph as a subset of \(X\times X\).

A subset of a standard Borel space is \emph{analytic} if it is the
projection of a Borel subset of its product with a Polish space.
A set is \emph{coanalytic} if its complement is analytic. We use
Souslin's theorem: a set which is both analytic and coanalytic is
Borel \cite{Kechris}.

\subsection{Reductions, universality, and smoothness}
For equivalence relations \(E\) on \(X\) and \(F\) on \(Y\), a
\emph{Borel reduction} from \(E\) to \(F\) is a Borel map
\(f:X\to Y\) satisfying
\[
 x\,E\,x'\quad\Longleftrightarrow\quad f(x)\,F\,f(x')
 \qquad(x,x'\in X).
\]
We write \(E\Bred F\) when such a map exists. Reductions compose.
The relation \({=}_Y\) is equality on \(Y\). An equivalence relation
is \emph{smooth} if it Borel reduces to \({=}_{2^{\N}}\), equivalently
if it admits complete invariants in a standard Borel space.
The relation \(E_0\) on \(2^{\N}\) is defined by
\[
 x\,E_0\,y\quad\Longleftrightarrow\quad
 \exists n\ \forall m\geq n\quad x(m)=y(m).
\]
The Harrington--Kechris--Louveau dichotomy says that every nonsmooth
Borel equivalence relation admits a Borel reduction of \(E_0\)
\cite{HKL}. We use it only for the classical low-level positive cases.

In the benchmark notation of \cite{Gao}, the symbol \(\ell_\infty\)
also denotes the equivalence relation on \(\R^{\N}\) given by
\[
 x\mathrel{\ell_\infty}y
 \quad\Longleftrightarrow\quad
 \sup_{n\in\N}|x(n)-y(n)|<\infty.
\]
Here \(\R^{\N}\) has its product topology. This is the relation
whose degree was mentioned in the introduction.

A \emph{common Borel upper bound} for a family \(\mathcal E\) is a
Borel equivalence relation \(F\) such that \(E\Bred F\) for every
\(E\in\mathcal E\). A universal \(\Gamma\) equivalence relation is
such an upper bound for all \(\Gamma\) equivalence relations, with
the additional requirement \(F\in\Gamma\). Unless a restriction is
stated, both definitions allow arbitrary Polish domains. No uniform
choice of the individual reductions is required.

\subsection{Baire category}
A set is \emph{nowhere dense} if its closure has empty interior,
\emph{meager} if it is a countable union of nowhere dense sets, and
\emph{comeager} if its complement is meager. The Baire category
theorem implies that a nonempty Polish space is not meager in itself.
Translations in a Polish group, being homeomorphisms, preserve
meager and comeager sets.

We use three standard category facts \cite{Kechris,Gao}. First, a
Borel map between Polish spaces is continuous on some comeager
\(G_\delta\) subspace. Second, the Kuratowski--Ulam theorem implies
that a Borel set \(D\subseteq X\times Y\) is comeager if and only if
\(D_x=\{y:(x,y)\in D\}\) is comeager for comeager many \(x\in X\).
Third, if \(P\) is standard Borel, \(X\) is Polish, and
\(D\subseteq P\times X\) is Borel, then
\[
 \{p\in P:D_p\text{ is comeager in }X\}
\]
is Borel. The last statement is Borel category quantification.
In the transfer argument it makes the conditions on a code Borel,
even though the existence of a suitable code is only analytic.

\subsection{Trees, summable matrices, and cosets}
Let \(S=\N^{<\N}\) be the set of finite sequences of natural numbers,
including the empty sequence \(o\). Write \(s^\frown n\) for the
sequence obtained by adjoining \(n\) to \(s\). A \emph{tree} is a
subset \(T\subseteq S\) closed under initial segments; the empty
tree is allowed. The set \(\Tr\) of trees is a closed subspace of
\(2^S\), hence compact. A branch of \(T\) is an \(x\in\N^{\N}\)
such that \(x{\upharpoonright}n\in T\) for every \(n\).
A tree is \emph{ill-founded} if it has a branch and
\emph{well-founded} otherwise. We denote these sets by \(\IF\) and
\(\WF\). The set \(\IF\) is analytic and non-Borel
\cite{Kechris,Marker}; it provides the contradiction in the proof.

For a countable set \(I\), the real Banach space \(\ell_1(I)\)
consists of the functions \(a:I\to\R\) for which
\(\|a\|_1=\sum_{i\in I}|a(i)|<\infty\). We use
\(A=\ell_1(S)\) and \(Z=\ell_1(S\times S)\), regarding elements
of \(Z\) as absolutely summable matrices. For \(a,b\in A\), set
\((a\otimes b)(s,t)=a(s)b(t)\); then
\(\|a\otimes b\|_1=\|a\|_1\|b\|_1\). Transposition is an
isometry of \(Z\). Only finite sums of these elementary tensors
and bounded coordinatewise operator actions will be needed.

A \emph{Polish group} is a topological group whose topology is Polish.
For a subgroup \(H\leq G\), our coset convention is
\[
 g\,E_H^G\,g'\quad\Longleftrightarrow\quad g^{-1}g'\in H.
\]
Thus the equivalence classes are the sets \(gH\).
No normality is assumed. If \(H\) is \(F_\sigma\), so is
\(E_H^G\), by continuity of \((g,g')\mapsto g^{-1}g'\).
For additive subgroups of Banach spaces, this is the relation
\(x-y\in H\); the sign is immaterial. In particular every real
linear subspace defines such an equivalence relation.

\section{A category transfer lemma}\label{sec:transfer}
The following lemma turns the existence of a common Borel upper
bound into regularity of a parameterized subgroup-membership test.
\begin{lemma}\label{lem:transfer}
Let \(P\) be a standard Borel space, let \(G\) and \(A_n\), \(n\in\N\),
be fixed Polish groups, and suppose
\[
 h_n:P\times A_n\longrightarrow G
\]
is Borel, with \(h_n(p,-)\) a continuous homomorphism for every \(p,n\).
Let \(H_p\) be the subgroup algebraically generated by these images.
For a Borel map \(z:P\to G\), put
\[
 M=\{p:z(p)\in H_p\}.
\]
If the relations \(E_{H_p}^G\), for \(p\notin M\), have a common
Borel equivalence-relation upper bound, then \(M\) is Borel.
\end{lemma}

\begin{proof}
First, \(M\) is analytic: membership has a finite-word witness in the
images of the countably many homomorphisms.

Fix a Borel equivalence relation \(F\) on a nonempty Polish space \(Y\).
There is a Polish code space \(C\) and a jointly Borel family
\(f_c:G\to Y\) representing every Borel map modulo a meager set.
Here is an explicit construction. Fix a countable basis of \(G\), a
dense sequence \((y_j)\) in \(Y\), and a complete compatible metric on
\(Y\). A code specifies an array of open sets \((U_{n,j})_{n,j\in\N}\).
Let \(s_{c,n}(g)\) be \(y_j\) for the least \(j\) with \(g\in U_{n,j}\),
or \(y_0\) if there is none. Define \(f_c(g)\) as the limit of these
simple maps when the sequence is Cauchy, and as \(y_0\) otherwise.
This evaluation is jointly Borel.

To represent a Borel \(f\), take a comeager \(G_\delta\) set \(D\)
on which \(f\) is continuous. Choose ambient open \(U_{n,j}\) whose
intersection with \(D\) is
\(D\cap f^{-1}(B(y_j,2^{-n}))\). On \(D\) the first eligible index
always gives an approximation within \(2^{-n}\), so \(f_c=f\) there.

Call \((p,c)\) a witness if
\begin{align}
 &\{(g,a): f_c(g)\,F\,f_c(g h_n(p,a))\}
       \text{ is comeager in }G\times A_n
       &&\text{for every }n, \label{eq:genhom}\\
 &\{g:\neg(f_c(g)\,F\,f_c(gz(p)))\}
       \text{ is comeager in }G. \label{eq:sep}
\end{align}
Both conditions are Borel in \((p,c)\), by Borel category
quantification. Hence the set \(W_F\) of parameters admitting a
witness is analytic.

We claim \(W_F\cap M=\varnothing\). Fix \(p,c\) satisfying
\eqref{eq:genhom}, fix \(n\), and abbreviate \(h_n(p,-)\) by \(h\).
For any fixed \(a_0\in A_n\), the map
\[
 (g,a)\longmapsto (g h(a_0),a_0^{-1}a)
\]
is a homeomorphism of \(G\times A_n\). On a comeager set, both
\[
 f_c(g)\,F\,f_c(g h(a)),\qquad
 f_c(g h(a_0))\,F\,f_c(g h(a_0)h(a_0^{-1}a))
                 =f_c(g h(a))
\]
hold. Transitivity and Kuratowski--Ulam give
\[
 \{g:f_c(g)\,F\,f_c(g h(a_0))\}\text{ comeager in }G.
\]
Intersecting finitely many translates of such sets shows the same
conclusion for every fixed element of \(H_p\). If \(z(p)\in H_p\),
this contradicts \eqref{eq:sep}.

Now suppose \(F\) bounds every off-\(M\) relation. For \(p\notin M\),
choose a Borel reduction \(f:E_{H_p}^G\to F\), and choose \(c\) such
that \(f_c=f\) on a comeager Borel set \(D\). For each fixed \(a\),
both \(g\) and \(g h_n(p,a)\) lie in \(D\) for comeager \(g\).
The reduction property and Kuratowski--Ulam give \eqref{eq:genhom}.
Since \(z(p)\notin H_p\), the reduction property on \(D\cap Dz(p)^{-1}\)
gives \eqref{eq:sep}. Thus \(P\setminus M\subseteq W_F\).
We have proved \(W_F=P\setminus M\). Both \(M\) and its complement are
analytic, so \(M\) is Borel.
\end{proof}

The lemma does not require a Borel choice of the reductions to \(F\).
It also imposes no Borelness assumption on \(H_p\) when \(p\in M\).

\section{An \(F_\sigma\) symmetric-tensor image}\label{sec:tensor}
Let \(S\) be a countably infinite set with a designated element \(o\),
and let
\[
 A=\ell_1(S;\R),\qquad Z=\ell_1(S\times S;\R).
\]
For \(a,b\in A\), write \(a\otimes b\) for the matrix
\((a(s)b(t))_{s,t}\). Let \(\mathcal S\subset Z\) be the real vector
space of symmetric finite-rank matrices. Finite rank here means
scalar matrix rank. It is equivalent to being an algebraic finite
sum of tensors of \(\ell_1\) vectors: a nonzero coordinate minor gives
the usual selected-row and selected-column factorization.
The condition \(\rank R\leq r\) is closed, since it is the vanishing
of every \((r+1)\)-minor.
We write \(\sym(a,b)=(a\otimes b+b\otimes a)/2\).

For a bounded linear \(D:A\to A\), its coordinatewise tensor actions
\(D\otimes\id\) and \(\id\otimes D\) are bounded on \(Z\).
Put
\[
 \mathcal C_D
 =\{(K(o,o),(D\otimes\id)K):K\in\mathcal S\}
 \subseteq \R\times Z.
\]

\begin{lemma}\label{lem:tensor}
If \(D\) is injective, then \(\mathcal C_D\) is an \(F_\sigma\)
real linear subspace of \(\R\times Z\).
\end{lemma}

\begin{proof}
For finite-rank \(R\), put \(Q=(\id\otimes D)R\). We first show
\[
 R=(D\otimes\id)K\text{ for some }K\in\mathcal S
 \quad\Longleftrightarrow\quad Q=Q^{\mathsf t}. \tag{*}
\]
The forward implication is immediate. For the reverse, choose a
minimal representation \(R=\sum_{i=1}^r u_i\otimes a_i\), with both
families linearly independent. Injectivity of \(D\) preserves the
independence of the \(Da_i\). Thus the first factor space of \(Q\)
is \(\operatorname{span}(u_i)\), and its second factor space is
\(D(\operatorname{span}(a_i))\). Symmetry equates them, so each \(u_i\)
has a preimage in \(\operatorname{span}(a_i)\). This gives an algebraic
lift \(K\). The symmetry of \((D\otimes D)K\), and injectivity on
algebraic tensors, imply that \(K\) is symmetric. The lift is unique.

We give explicit closed charts for its scalar coordinate. Fix
\(r\geq1\), distinct indices \(s_1,\ldots,s_r\), distinct indices
\(q_1,\ldots,q_r\), and \(k\geq1\). Define
\[
 a_i(R)=R(s_i,-),\qquad
 N(R)_{ji}=(D a_i(R))(q_j).
\]
On \(\det N(R)\ne0\), set
\[
 b_i(R)=\sum_{j=1}^r (N(R)^{-1})_{ij}R(q_j,-),\qquad
 \kappa(R)=\sum_{i=1}^r a_i(R)(o)b_i(R)(o).
\]
Consider the set of \((\lambda,R)\) satisfying
\[
 \rank R\leq r,\qquad Q=Q^{\mathsf t},\qquad
 |\det N(R)|\geq 1/k,\qquad \lambda=\kappa(R). \tag{**}
\]
This set is closed. All row maps and determinants are continuous,
and \(\kappa\) is continuous on the open nonzero-determinant domain;
the displayed determinant bound preserves its graph under limits.

The selected rows are independent on this chart. They form a basis
of the row space, and (*) says that the column space is their image
under \(D\). Solving at the coordinates \(q_j\) gives
\[
 R=\sum_{i=1}^r Da_i(R)\otimes b_i(R).
\]
Consequently \(K=\sum_i a_i(R)\otimes b_i(R)\) is its unique
symmetric lift, with \(K(o,o)=\kappa(R)\). Thus (**) is contained
in \(\mathcal C_D\).

Conversely, every nonzero lifted matrix admits independent selected
rows, a nonzero coordinate minor of their \(D\)-images, and a positive
integer determinant bound. It belongs to one of these countably
many charts. When \(R=0\), injectivity forces \(K=0\). Hence
\(\mathcal C_D\) is the union of these charts and the singleton
\(\{(0,0)\}\).
\end{proof}

The lemma does not assume that \(D\) has closed range or a bounded
inverse. Symmetry forces range membership inside a finite-dimensional
space that can be recovered from the rows of \(R\).

\section{Tree kernels and the fixed Polish group}\label{sec:trees}
Now take \(S=\N^{<\N}\), with root \(o=\varnothing\). Let \(\Tr\)
be the compact space of downward-closed trees on \(S\). For \(T\in\Tr\),
define
\[
 (S_Ta)(s)=\sum_{\{n:s^\frown n\in T\}}a(s^\frown n),
 \qquad D_T=I-2S_T.
\]
Absolute summability gives \(\|S_T\|\leq1\) and \(\|D_T\|\leq3\).
The map \((T,a)\mapsto D_Ta\) is jointly continuous: for fixed
finitely supported \(a\), only finitely many bits of \(T\) are used,
and the common norm bound controls the remaining tail.

If \(x\) is a branch of \(T\), the vector
\[
 a(x{\upharpoonright}n)=2^{-n},\qquad a(s)=0
       \text{ off this branch}
\]
has norm two, satisfies \(D_Ta=0\), and has \(a(o)=1\).
Conversely, a nonzero coordinate of a kernel vector satisfies
\[
 a(s)=2\sum_{\{n:s^\frown n\in T\}}a(s^\frown n),
\]
so it has a nonzero child in \(T\). Iterating yields a branch.
Therefore
\[
 T\in\IF
 \ \Longleftrightarrow\ \ker D_T\ne0
 \ \Longleftrightarrow\
 \exists a\in\ker D_T\ (a(o)=1). \tag{1}
\]
In particular \(D_T\) is injective when \(T\in\WF\).

Let \(X=A\times A\) and \(C=\R\times Z\). Define the continuous
bilinear map
\[
 \beta((a,u),(b,v))
 =\left(a(o)b(o),\frac{u\otimes b+v\otimes a}{2}\right)
 \in C.
\]
On \(G=X\times X\times C\), define
\[
 (x,y,c)(x',y',c')
 =(x+x',y+y',c+c'+\beta(x,y')).
 \tag{2}
\]
Bilinearity gives associativity, and
\[
 (x,y,c)^{-1}=(-x,-y,-c+\beta(x,y)).
\]
Thus \(G\) is one fixed two-step nilpotent Polish group.
Here two-step nilpotence means that every commutator is central;
in this group all commutators lie in \(\{0\}\times\{0\}\times C\).

For each \(T\), consider the two homomorphisms from \((A,+)\):
\[
 h_{T,0}(a)=((a,D_Ta),0,0),\qquad
 h_{T,1}(b)=(0,(b,D_Tb),0).
 \tag{3}
\]
Their images are closed graph subgroups, and both families are
jointly continuous. Let \(H_T\) be the subgroup they generate.
The mixed commutator is the central element with value
\[
 \beta((a,D_Ta),(b,D_Tb))
 =\left(a(o)b(o),(D_T\otimes\id)
           \frac{a\otimes b+b\otimes a}{2}\right).
 \tag{4}
\]
Every symmetric finite-rank tensor is a finite real linear
combination of the displayed symmetric tensors. Real multiples
are realized by scaling one commutator input. Therefore the
commutators generate exactly \(\mathcal C_{D_T}\), and
\[
 H_T=\graph(D_T)\times\graph(D_T)\times\mathcal C_{D_T}
 \quad\text{as a set of coordinates in }G. \tag{5}
\]
Indeed the right-hand side is a subgroup containing both axes;
conversely, axis elements give the first two coordinates and
commutators correct the central coordinate.

For \(T\in\WF\), Lemma~\ref{lem:tensor} and (5) show that \(H_T\)
is \(F_\sigma\). Hence \(E_{H_T}^G\) is \(F_\sigma\), since group
difference is continuous.

Fix the central element
\[
 z=(0,0,(1,0))\in G.
\]
By (1), a branch vector \(a\) gives the symmetric tensor \(a\otimes a\)
with central image \((1,0)\). Conversely, if
\[
 K(o,o)=1,\qquad (D_T\otimes\id)K=0,
\]
then the \(o\)-th column of \(K\) is a kernel vector with root
coordinate one. We have proved the exact equivalence
\[
 z\in H_T\quad\Longleftrightarrow\quad T\in\IF. \tag{6}
\]

\begin{proof}[Proof of Theorem~\ref{thm:main}]
Suppose a Borel equivalence relation \(F\) bounded every \(F_\sigma\)
equivalence relation. It would bound \(E_{H_T}^G\) for every
well-founded \(T\). Apply Lemma~\ref{lem:transfer} to the two jointly
continuous families (3), the fixed group \(G\), and the constant
element \(z\). Equation (6) would make \(\IF\) Borel, contradicting
its classical analytic non-Borelness.
\end{proof}

\section{Linear coset relations on a fixed Banach space}\label{sec:linear}
The noncommutative group above organizes the bilinear terms. Its
coset relations can nevertheless be represented by additive cosets
of real linear subspaces, as the next explicit reduction shows.

\begin{corollary}\label{cor:linear}
Even the coset equivalence relations of \(F_\sigma\) real linear
subspaces of one fixed separable Banach space isomorphic to \(\ell_1\)
have no common Borel equivalence-relation upper bound.
\end{corollary}

\begin{proof}
For fixed \(D\), write \(g=((a,u),(b,v),c)\in G\), and put
\[
 r=u-Da,\quad s=v-Db,\quad
 \theta_D(g)=\left(r,s,c-\left(0,\frac{r\otimes b}{2}\right)\right).
\]
This continuous map reduces the coset relation for
\(\graph(D)\times\graph(D)\times\mathcal C_D\) to the additive
coset relation for \(\{0\}^2\times\mathcal C_D\) on \(A^2\times C\).
Indeed, membership of \(g^{-1}g'\) in that subgroup first requires
\(r=r'\) and \(s=s'\). Writing \(e=b'-b\), the remaining central
condition is
\[
 c'-c-\beta((a,u),(e,De))\in\mathcal C_D.
\]
Since
\[
 \beta((a,u),(e,De))
 =\left(a(o)e(o),(D\otimes\id)\sym(a,e)\right)
       +\left(0,\frac{r\otimes e}{2}\right),
\]
this is precisely equality of the last coordinates of
\(\theta_D(g)\) and \(\theta_D(g')\) modulo \(\mathcal C_D\).
Conversely the continuous map
\((r,s,\eta)\mapsto((0,r),(0,s),\eta)\) is a reduction in the other
direction.

For well-founded \(T\), the subspace
\(\{0\}^2\times\mathcal C_{D_T}\) is \(F_\sigma\). A common Borel
bound for all these linear cosets would therefore bound all the
relations used in Theorem~\ref{thm:main}, a contradiction.
The fixed Banach space
\(A^2\times C=\ell_1(S)^2\times\R\times\ell_1(S\times S)\),
with its finite sum norm, is linearly isomorphic to \(\ell_1\).
\end{proof}

\section{Universality in the Borel hierarchy}\label{sec:hierarchy}
\subsection{Arbitrary Polish spaces}
\begin{corollary}\label{cor:hierarchy}
For \(1\leq\alpha<\omega_1\), the existence of a universal equivalence
relation on arbitrary Polish spaces is as follows:
\[
\begin{array}{c|c|c}
 \alpha & \Sigma^0_\alpha & \Pi^0_\alpha\\ \hline
 1 & \text{yes: equality on }\N
   & \text{yes: equality on }2^\N\\
 2 & \text{no}
   & \text{yes: equality on }2^\N\\
 \alpha\geq3 & \text{no} & \text{no}.
\end{array}
\]
\end{corollary}

\begin{proof}
Every \(\Sigma^0_\alpha\) class with \(\alpha\geq2\), and every
\(\Pi^0_\alpha\) class with \(\alpha\geq3\), contains all
\(F_\sigma\) sets. A universal relation in any of these classes
would be a Borel upper bound forbidden by Theorem~\ref{thm:main}.

An open equivalence relation has open classes, hence only countably
many classes in a second-countable space. Enumerating them gives a
Borel reduction to equality on the discrete space \(\N\), which is
itself open.

Closed and \(G_\delta\) equivalence relations are smooth, so they
reduce to equality on \(2^\N\). For completeness this follows from
the Harrington--Kechris--Louveau dichotomy \cite{HKL}. If a
\(G_\delta\) equivalence relation were nonsmooth, there would be a
Borel reduction of \(E_0\) to it. Restrict the reduction to an
\(E_0\)-invariant comeager \(G_\delta\) set on which it is continuous:
intersect a continuity domain with all its finite-coordinate-flip
translates. The pullback relation is \(G_\delta\) on this domain.
Each of its \(E_0\)-classes is countable and dense there, so, being
\(G_\delta\), is comeager there. Two distinct such classes would be
disjoint comeager sets in a nonempty Polish space, a contradiction.
Equality on \(2^\N\) is closed and hence also \(G_\delta\).
\end{proof}

The usual hierarchy begins at level one. If level zero is defined
to consist of the clopen sets on both sides, equality on discrete
\(\N\) is universal at that level as well.

This conclusion concerns arbitrary Polish domains. It does not
contradict Rosendal's universal \(K_\sigma\) equivalence relation
\cite{Rosendal}; \(K_\sigma\) is a more restrictive requirement than
\(F_\sigma\) on a noncompact Polish space.

\subsection{The Cantor-space formulation}
To relate the preceding conclusions to the formulation in
\cite[Section~9]{KanoveiLyubetsky}, we record the small increase in
Borel complexity needed to move an arbitrary Polish domain into
Cantor space.

\begin{proposition}\label{prop:cantor-extension}
Every \(F_\sigma\) equivalence relation on a Polish space is Borel
reducible to a \(\Delta^0_3\) equivalence relation on \(2^{\N}\).
\end{proposition}

\begin{proof}
Let \(E\) be \(F_\sigma\) on \(X\). The assertion is immediate
when \(X\) is empty, so assume otherwise. Equip \(X\) with a finer
zero-dimensional Polish topology having the same Borel sets, and
embed this space homeomorphically as a \(G_\delta\) subset
\(G\subseteq2^{\N}\). These are standard facts about Polish spaces
\cite{Kechris}. The identity to the original topology is continuous,
so \(E\) remains \(F_\sigma\) in the finer topology. Its transported
relation \(\widetilde E\) is therefore \(F_\sigma\) relative to
\(G^2\). Choose an ambient \(F_\sigma\) set
\(B\subseteq(2^{\N})^2\) such that
\(\widetilde E=B\cap G^2\), and put
\[
 E^*=(B\cap G^2)\cup\Delta_{2^{\N}},
 \qquad
 \Delta_{2^{\N}}=\{(x,x):x\in2^{\N}\}.
\]
This is an equivalence relation: it agrees with \(\widetilde E\)
on \(G\) and has singleton classes outside \(G\). The embedding,
viewed as a Borel map from the original \(X\), reduces \(E\) to
\(E^*\).

Both an \(F_\sigma\) set and a \(G_\delta\) set belong to
\(\Delta^0_3\). The class \(\Delta^0_3\) is closed under finite
unions and intersections, and the diagonal is closed. Hence
\(E^*\in\Delta^0_3((2^{\N})^2)\).
\end{proof}

\begin{corollary}\label{cor:cantor}
The \(\Delta^0_3\) equivalence relations on Cantor space have no
common Borel equivalence-relation upper bound. In particular, for
every countable \(\alpha\geq3\), neither the \(\Sigma^0_\alpha\)
nor the \(\Pi^0_\alpha\) equivalence relations on Cantor space have
a universal member.
\end{corollary}

\begin{proof}
A common Borel bound for these \(\Delta^0_3\) relations would, by
Proposition~\ref{prop:cantor-extension} and transitivity of Borel
reducibility, bound every \(F_\sigma\) equivalence relation on every
Polish space. This contradicts Theorem~\ref{thm:main}. For
\(\alpha\geq3\), each of the two indicated pointclasses contains
\(\Delta^0_3\), giving the last assertion.
\end{proof}

Thus the higher-level negative answer holds in the Cantor-space
setting discussed by Kanovei and Lyubetsky, as well as on arbitrary
Polish domains. At level \(\Sigma^0_2\), however, every equivalence
relation on a compact Polish space is \(K_\sigma\) whenever it is
\(F_\sigma\). Rosendal's compact representative of the universal
\(K_\sigma\) degree \cite[Definition~18 and Proposition~19]{Rosendal}
therefore supplies a universal member in that compact setting.
The intersection with \(G^2\) in Proposition~\ref{prop:cantor-extension}
need not preserve \(F_\sigma\), which is precisely why the argument
does not transfer the arbitrary-domain level-two conclusion to
compact domains.

\noindent\begin{minipage}{\textwidth}
\small
\textsc{School of Mathematical Sciences and School of Pre-university,}\\
\textsc{Dalian Minzu University}\\[0.4em]
\textit{Email address:}\enspace
\href{mailto:wangruijun@dlnu.edu.cn}{\texttt{wangruijun@dlnu.edu.cn}}
\end{minipage}
\end{document}